\documentclass[10pt, conference, letterpaper]{IEEEtran}
\usepackage{amsmath}
\usepackage{amssymb}
\usepackage{amsthm}
\usepackage{bm}
\usepackage{dsfont}
\usepackage{subcaption}
\usepackage[ruled,vlined,linesnumbered]{algorithm2e}  
\usepackage{url}
\usepackage{xcolor}
\usepackage{mdframed}
\usepackage{graphicx}

\begin{document}
\newcommand{\vect}[1]{\mathbf{{#1}}}
\newcommand{\set}[1]{\mathcal{{#1}}}
\newcommand{\abs}[1]{\left|{#1}\right|}
\newcommand{\ceiling}[1]{\lceil{#1}\rceil}
\newcommand{\floor}[1]{\lfloor{#1}\rfloor}

\newcommand{\prtc}[1]{\left\{{#1}\right\}}
\newcommand{\prts}[1]{\left[{#1}\right]}
\newcommand{\prtr}[1]{\left({#1}\right)}
\newcommand{\state}[1]{\mathcal{S}({#1})}
\newcommand{\schedu}[1]{\mathcal{X}({#1})}
\newcommand{\stateopt}[1]{\mathcal{S}^*({#1})}
\newcommand{\scheduopt}[1]{\mathcal{X}^*({#1})}

\newcommand{\policy}{\pi}

\newcommand{\iterstep}{k}

\newcommand{\pidle}[1]{R({#1})}

\newcommand{\avgwheeltime}{root-mean-square wheel time}
\newcommand{\avgwt}{\bar{\Omega}}

\newcommand{\productwheel}{product wheel}
\newcommand{\wheeltime}{wheel time}
\newcommand{\wt}{\Omega}
\newcommand{\wttime}[1]{\Omega_{#1}}
\newcommand{\numitem}{N}
\newcommand{\pitem}[1]{I_{#1}}

\newcommand{\batch}[1]{batch}
\newcommand{\optimalbatch}{\bm{\lambda}}
\newcommand{\vecbatch}{\bm{\lambda}}
\newcommand{\optimalbatcht}[1]{\bm{\lambda}(#1)}
\newcommand{\optimalbatchtp}[1]{\bm{\lambda'}(#1)}
\newcommand{\fwheeltime}[1]{\bm{J}(#1)}
\newcommand{\bat}[1]{\lambda_{#1}}
\newcommand{\statesa}[1]{\avgwt(#1)}
\newcommand{\statesap}[1]{\avgwt'(#1)}
\newcommand{\coolingconst}{C}

\newcommand{\tempsa}[1]{T(#1)}
\newcommand{\tmaxsa}{W}
\newcommand{\feasiblew}{feasible wheels}
\newcommand{\feasiblewset}{\mathcal{S}}
\newcommand{\proposf}{f}

\newcommand{\unitbat}[1]{M_{#1}}

\newcommand{\timebat}[1]{t_{#1}}

\newcommand{\numslot}{H}
\newcommand{\timeslot}{T}

\newcommand{\setnode}{\mathcal{N}}
\newcommand{\setslot}{\mathcal{H}}

\newcommand{\proditem}[2]{s_{#1}^{#2}}

\newcommand{\inventoryc}{inventory cost}
\newcommand{\changeoverc}{set-up cost}
\newcommand{\triggerpoint}{trigger point}
\newcommand{\costtolerance}{cost tolerance}
\newcommand{\costconstraint}{cost constraint}
\newcommand{\inventoryconstraint}{inventory constraint}
\newcommand{\indexset}{\mathcal{N}}
\newcommand{\timeset}{\mathcal{H}}

\newcommand{\ct}{\tau}

\newcommand{\tp}[1]{Q_{#1}}
\newcommand{\cost}[2]{c_{#1}^{#2}}
\newcommand{\icost}[1]{k_{#1}}
\newcommand{\ccost}[1]{p_{#1}}

\newcommand{\demand}[2]{D_{#1}^{#2}}
\newcommand{\alldemand}{\mathcal{D}}
\newcommand{\invent}[2]{I_{#1}^{#2}}

\newcommand{\backoffleftname}{back-off time left}
\newcommand{\backoffleft}[2]{L_{{#1}}({{#2}}) }

\newcommand{\stategreedyname}{state}
\newcommand{\stategreedy}[1]{St({{#1}}) }

\newcommand{\numtransmit}[1]{N_{{#1}}}

\newcommand{\tx}[2]{x_{{#1}}({#2})}
\newcommand{\gentime}[2]{U_{{#1}}({#2})}
\newcommand{\procage}[2]{A_{{#1}}({#2})}

\newcommand{\weight}[2]{W_{{#1}}({#2})}

\newcommand{\tiematprocagemax}[1]{t_{{#1}}^*}

\newcommand{\figref}[1]{Figure~\ref{#1}}

\newcommand\bluesout{\bgroup\markoverwith{\textcolor{blue}{\rule[0.5ex]{2pt}{0.4pt}}}\ULon}
\newcommand\magentasout{\bgroup\markoverwith{\textcolor{green}{\rule[0.5ex]{2pt}{0.4pt}}}\ULon}
\newcommand\redsout{\bgroup\markoverwith{\textcolor{red}{\rule[0.5ex]{2pt}{0.4pt}}}\ULon}
\newcommand\purplesout{\bgroup\markoverwith{\textcolor{purple}{\rule[0.5ex]{2pt}{0.4pt}}}\ULon}

\newtheorem{theorem}{Theorem}
\newtheorem{corollary}{Corollary}
\newtheorem{lemma}{Lemma}
\newtheorem{property}{Property}
\newtheorem{definition}{Definition}

\newtheorem*{proper}{\algname}
\surroundwithmdframed{proper}

\newtheorem*{power}{\algnamepower}
\surroundwithmdframed{power}

\newtheorem*{powerc}{\algnamepowerc}
\surroundwithmdframed{powerc}
% \newmdtheoremenv{proper}{\algname}

%%% Comment
\newcommand{\rob}[1]{\todo[author=Rob,color=blue!30,inline]{#1}}
\newcommand{\wasin}[1]{\todo[author=WS,color=magenta!30,inline]{#1}}

%%% Add text
\newcommand{\robadd}[1]{{\color{blue}{#1}}}
\newcommand{\wasinadd}[1]{{\color{magenta}{#1}}}

	%\title{Optimizing Age of Information with Processing Constraints}
	\title{Optimizing Product Wheel Time in Lean Manufacturing Systems}

% 	\author{\IEEEauthorblockN{1\textsuperscript{st} Wasin Meesena}
% \IEEEauthorblockA{\textit{Department of Mathematics and Statistics} \\
% \textit{Carleton College}\\
% MN, United States \\
% meesenat@carleton.edu}
% \and 
% \IEEEauthorblockN{2\textsuperscript{nd} Robert Thompson}
% \IEEEauthorblockA{\textit{Department of Mathematics and Statistics} \\
% \textit{Carleton College}\\
% MN, United States \\
% rthompson@carleton.edu}}
\author{\IEEEauthorblockN{Wasin~Meesena\IEEEauthorrefmark{2}, 
                          Robert~Thompson\IEEEauthorrefmark{2}}
\IEEEauthorblockA{\IEEEauthorrefmark{2}Department of Mathematics and Statistics, Carleton College, Northfield, Minnesota, USA\\
(meesenat@carleton.edu, rthompson@carleton.edu)}
}

	\maketitle\thispagestyle{empty}

    \begin{abstract}

Lean manufacturing is a production method focused on reducing production times, eliminating waste, and synchronizing production with fluctuating demand. A standard lean manufacturing methodology is the \textit{product wheel}, a repeating sequence of production of various items. If this product wheel sequence is short, it is easier to interrupt or alter production to adjust for failures or fluctuations in demand, so the manufacturing process is leaner.  However, a sequence that is too short results in frequent changeover from the production of one item to the next, yielding higher costs. This study formulates the product wheel methodology as an optimization problem and proposes two approaches to solving this problem: one via a relaxation to integer linear programming, and another via the probabilistic optimization technique of simulated annealing. We assess the performance of these two approaches through simulations and analyze the  relationships between production leanness and costs.
\end{abstract}


	\IEEEpeerreviewmaketitle

    \section{Introduction}

Economies of scale are foundational to traditional manufacturing, \cite{tradition}.   High-volume operations help minimize total unit costs, but result in large inventories and other inefficiencies.  By contrast, lean manufacturing -- pioneered by Toyota in the 1950s --  focuses on minimizing waste and being highly responsive to demand, \cite{gupta2013literature, lean}. By producing items in harmony with the demand rate, lean manufacturing reduces overproduction, waiting times, excess inventory, and the impacts of defects and demand fluctuations,  \cite{king2013product,womack2007machine}.  Many businesses, particularly those in the automotive and aerospace industries, embrace lean manufacturing, \cite{day2004learning, jones1999seeing, womack1999lean}.

The \textit{product wheel} is a concept derived from the theory of lean manufacturing, \cite{wheel-lean,boonkanok2021consumer}.  The core idea of the product wheel is to repeat a sequence of production of various items, with a fixed quantity of each item produced at each step in the sequence.  This methodology has two key advantages: (i) the fixed nature of the production sequence simplifies operations, and (ii) production can more easily be interrupted to address defects or demand fluctuations. Furthermore, the \textit{shorter} the production sequence, the easier it is to interrupt or modify production, resulting in a \textit{leaner} manufacturing process.  However, a shorter sequence  leads to more frequent changeover from the production of one item to the next, incurring greater costs, \cite{king2013product}.  The present work focuses on developing an optimization model to quantitatively determine the trade-off between leanness and production cost in the product wheel methodology.

The contributions of this paper are as follows.
\begin{itemize}
     \item As far as we know, this is the first work that formulates and aims to optimize leanness, which is  the time required for completing a product wheel cycle, in the product wheel methodology under demand and cost constraints.

     \item  Two approaches to the product wheel optimization problem are described: (i) the probabilistic optimization technique of simulated annealing for optimization in the general nonconvex case where skipping products in the product wheel sequence is allowed, and (ii) an integer linear programming approach for a relaxed linear version of the problem, where product skipping is prohibited and the product wheel stays fixed through the entire manufacturing process.

     \item The performance of these two approaches is evaluated through simulations.  In particular, we investigate the relationships between leanness and different types of costs incurred.
    \end{itemize}

This paper is organized as follows. The problem formulation is presented in Section \ref{sec:formulation}. 
Section \ref{sec:approach} develops two approaches to solving this optimization problem. In Section \ref{sec:num}, these two approaches are tested via numerical simulation.  In Section \ref{sec:conclusion} we conclude by discussing further directions of research.
     
    \section{Problem Formulation}
\label{sec:formulation}

The formulation of our model is guided by three characteristics of an ideal product wheel: (i) the resulting production should meet demand, (ii) the total cost should not exceed a given \textit{cost tolerance}, and (iii) the \textit{wheel time} (total time for one complete cycle of the product wheel) should be minimized.

Consider a product wheel for manufacturing $N$ items over $\numslot$  time periods (e.g., 12 months), each with length $\timeslot$ (e.g., 100 hours).  Items and time periods are indexed by the sets $ \indexset= \prtc{1, \dotsc, N}$  and $\setslot =\prtc{1, 2, \dotsc, \numslot}$, respectively.  In one cycle of a \productwheel, we produce $\bat{i}$ batches of item $i$.  We do not produce partial batches, \cite{king2019lean}, so $\bat{i}\in \mathbb{N}.$  Each batch of item $i$ contains $\unitbat{i}$ units of the item, and the time needed to produce a single batch of item $i$ is $\timebat{i}$. Therefore, during one cycle of the  \productwheel, a total time of $\bat{i}\timebat{i}$ is spent on the production of item $i$.  The total time $\wt$ required to execute a full cycle of the product wheel, called the \textit{\wheeltime}, can be written 
\begin{equation}\label{wheelt}
    \wt = \sum_{i\in \indexset}\bat{i}\timebat{i}.
\end{equation}
If we do not allow partial runs of the product wheel, there are exactly $\floor{\frac{\timeslot}{\wt}}$ product wheel cycles possible in a single time period.

We now introduce the idea of item skipping.  During a given time period, it may not be necessary to produce all items on the product wheel.  Instead of modifying the product wheel itself to account for this, we introduce the possibility of skipping the production of a given item during that time period.  The choice of skipping is encoded by a $\textit{skipping coefficient}$ $\proditem{i}{h}$, where $\proditem{i}{h} = 1$ if item $i$ is to be produced in time period $h$, and $\proditem{i}{h} = 0$ if the item is skipped in the given time period.

Our wheel time then varies across time periods, and the  \textit{\wheeltime \ in period $h$}, denoted by $\wttime{h}$, is written as 
\begin{equation}\label{wheeltimepartial}
    \wttime{h} = \sum_{i\in \indexset} \bat{i}\timebat{i}\proditem{i}{h}.
\end{equation}
We complete $\floor{\frac{\timeslot}{\wttime{h}}}$ full cycles of the product wheel in time period $h$.  Note that if no skipping is allowed, all $\proditem{i}{h}$ are $1$ and $\wttime{h} = \wt$ for all time periods $h$.

To achieve lean manufacturing, the wheel time should be small.  Because item skipping causes wheel time to vary across time periods, we set our optimization objective to be the minimization of \avgwheeltime \  $\avgwt$, defined as 
\begin{equation}\label{avgwheel}
    \avgwt = \sqrt{\sum_{h \in \setslot} \wttime{h}^2/\numslot}.
\end{equation}
Root-mean-square wheel time is preferred to the average \wheeltime \ because it penalizes outliers more heavily -- a time period with a lengthy \wheeltime \ is highly undesirable.

Denote by $\demand{i}{h}$ the demand for item $i$ in period $h$ and by $\invent{i}{h}$ the inventory of item  $i$ after production in time period $h$ but before satisfying the demand $\demand{i}{h}$.  For simplicity, we assume that all demand $\demand{i}{h}$ is known before production planning. Current inventory can be expressed in terms of previous inventory, previous demand, and current production,
\begin{equation}\label{eq:costinventory}
    \invent{i}{h}= \invent{i}{h-1}-\demand{i}{h-1} +\bat{i}\unitbat{i}\proditem{i}{h}\floor{\frac{\timeslot}{ \wttime{h}}}.
\end{equation}

We consider two sources of costs: \textit{\changeoverc} \ and \textit{\inventoryc}. The \changeoverc \ $\ccost{i}$ is the cost incurred transitioning to production of item $i$ from the production of any other item.  Since there are $\proditem{i}{h}\floor{\frac{\timeslot}{ \wttime{h}}}$ product wheel cycles in period $h$, the total \changeoverc \ due to production of item $i$ in period $h$ is  $\ccost{i} \proditem{i}{h}\floor{\frac{\timeslot}{ \wttime{h}}}$. Note that  if we do not produce a particular item, no \changeoverc s are incurred for that item.  The \inventoryc \ $\icost{i}$ is the cost incurred when storing one unit of item $i$ for a single time period.  At the beginning of period $h$, there will be $\invent{i}{h-1}-\demand{i}{h-1}$ units of item $i$ that must be kept in inventory, so the total \inventoryc \ for item $i$ in period $h$ is $\icost{i}(\invent{i}{h-1}-\demand{i}{h-1}).$  The  cost of production $\cost{i}{h}$ associated with item $i$ and period $h$ is the combination of total \changeoverc \ and total \inventoryc,
\begin{equation}\label{costper}
    \cost{i}{h} = \ccost{i} \proditem{i}{h}\floor{\frac{\timeslot}{ \wttime{h}}} +\icost{i}(\invent{i}{h-1}-\demand{i}{h-1}).
\end{equation}

Finally, we have two requirements for production.  First, production must meet the demand for all items in all time periods:
\begin{align}\label{eq:inventoryconstr}
    \invent{i}{h}  \geq \demand{i}{h}.
\end{align}
Second, the total cost of production cannot exceed a fixed \textit{\costtolerance} \  $\ct$:
\begin{align}\label{eq:costconstr}
    \sum_{h \in \setslot}\sum_{i\in \indexset} \cost{i}{h} \leq \ct.
\end{align}
To create a lean manufacturing plan, we want to minimize the \avgwheeltime\ subject to these constraints.  Thus we arrive at our central optimization problem:
\begin{subequations}
\label{opt:problem}
\begin{align}
   &\text{minimize}  && \ \avgwt   \nonumber\\
     &\text{s.t. } && \sum_{h \in \setslot}\sum_{i\in \indexset} \cost{i}{h} \leq \ct   \   \label{opt:constr1}\\
      & \ &&  \invent{i}{h}\geq \demand{i}{h}\ &&&   i \in \setnode,  h \in \setslot; \label{opt:constr2}\\
         & \  && \bat{i}\in \mathbb{N}, \proditem{i}{h} \in \{0,1\}  \ &&&   i \in \setnode,  h \in \setslot. \label{opt:constr3}
\end{align}
\end{subequations}
    \section{Optimization Approach}
\label{sec:approach}

In the optimization problem \eqref{opt:problem}, there are $2^{N\cdot H}$ choices for skipping coefficients $\prtc{\proditem{i}{h}}_{ i \in \setnode}^{  h \in \setslot}$, and for each of these choices (if feasible) there is an optimal product wheel to be computed.  Searching for the optimal wheel time over all possible skipping choices is intractable, so we consider two restricted choices of skipping coefficients.

First, we consider optimization with no item skipping, so $\proditem{i}{h} = 1$ for all items and time periods. Wheel time will then be the same for all time periods, and with a small relaxation, \eqref{opt:problem} can be reformulated as an integer linear program.  Second, we consider non-anticipative item skipping, where skipping coefficients in a given time period are determined by a criterion applied to the current production variables.  A heuristic optimization based on simulated annealing is then used to solve this simplified problem.

\subsection{No item skipping}

With no item skipping, all $\proditem{i}{h} = 1$.  The objective function in \eqref{opt:problem} then becomes identical to \eqref{wheelt}, and thus linear in the $\lambda_i$.  To arrive at a linear program, we relax equation \eqref{eq:costinventory} by allowing a fractional number $\frac{\timeslot}{\wt}$ of cycles of the product wheel in a given time period.  With this relaxation \eqref{eq:costinventory} becomes

\begin{equation}\label{eq:costinventory_relaxed}
    \invent{i}{h}= \invent{i}{0}- \sum_{j=1}^{h-1} \demand{i}{j}+  \bat{i}h \unitbat{i}\frac{ \timeslot}{\sum_{i\in \setnode} \bat{i}\timebat{i}}.
\end{equation}
After expressing inventory via \eqref{eq:costinventory_relaxed} and applying some algebraic manipulation and the relaxation assumption, the cost and demand constraints \eqref{opt:constr1} and \eqref{opt:constr2} are also seen to be linear.  The optimization problem $\eqref{opt:problem}$ can then be recast as an integer linear program and solved with standard software such as Python-MIP, \cite{mip}.

\subsection{Non-anticipative item skipping}

In practice, items can be skipped when the current inventory does not require their production.  For example, item production in time period $h$ may be skipped if remaining inventory meets the next period's demand, i.e. $\invent{i}{h-1}-\demand{i}{h-1} \geq  \demand{i}{h}$.  Item production could also be triggered when the remaining inventory falls below a set value $\tp{i}$ known as a \textit{trigger point}.  Combining these criterion results in a skipping coefficient

\begin{align}\label{eq:trigger}
     \proditem{i}{h}=\begin{cases} 
      0       & \text{if $\invent{i}{h-1}-\demand{i}{h-1} \geq \max{(\tp{i}, \demand{i}{h})}$,}\\
      1    & \text{otherwise.} \\
   \end{cases}
\end{align}
Skipping according to \eqref{eq:trigger} is \textit{non-anticipative}, relying on a criterion applied only to the current (or previous) production values, with no need for knowledge of future states, \cite{sun2019age}.  Applying a non-anticipative skipping policy simplifies the solution of the optimization problem \eqref{opt:problem} because all skipping coefficients can be determined by a choice of product wheel.

  We  outline in Algorithm \ref{alg:simannealing} the application of  simulated annealing to determine an optimal \productwheel \ for a non-anticipatory skipping policy.   The central idea behind simulated annealing is to iteratively search for a more optimal state with a guessing strategy that ``cools'', becoming more restrictive as the iteration progresses, \cite{simulatedannealref}.  We apply a constant cooling $\coolingconst$ over a fixed number $\tmaxsa$ of iterations.  Write $\optimalbatcht{\iterstep} = \begin{bmatrix} \bat{1}(k) & \bat{2}(k) & \dots & \bat{\numitem}(k) \end{bmatrix}$ and $\statesa{\iterstep}$ for the  \productwheel \ and \avgwheeltime \ values at iteration  $k$, and suppose that $\optimalbatcht{0}$ and the demand schedule $\alldemand$ are known.  Let $\pi$ be shorthand for our non-anticipative skipping policy.

\begin{algorithm}
\caption{Simulated Annealing}\label{alg:simannealing}
\SetKwInOut{Input}{Input}\SetKwInOut{Output}{Output}
\Input{$\alldemand, \coolingconst, \policy,\tmaxsa, \optimalbatcht{0} $}
\Output{An optimal batch $\optimalbatch^*$}
  \SetAlgoLined
\While{$ \iterstep<\tmaxsa$ }{
 Propose $\optimalbatchtp{\iterstep+1} $ given the current  $\optimalbatcht{\iterstep}$ \label{inalg:propose}\\
\While{ $\optimalbatchtp{\iterstep+1}$ \text{is not feasible} \label{inalg:check}}{
Propose $\optimalbatchtp{\iterstep+1} $ given the current $\optimalbatcht{\iterstep}$}\label{inalg:check2}
\If{$\statesap{\iterstep+1} \leq \statesa{\iterstep}$ \label{inalg:betterstate}}{$ \optimalbatcht{\iterstep+1} \gets  \optimalbatchtp{\iterstep+1}$} \label{inalg:betterstate2}
\Else{Sample $z$ from a uniform distribution on $[0,1]$ \label{inalg:worsestate}\\
\If{$z<\exp{[- \prtr{\statesap{\iterstep+1} -\statesa{\iterstep} }/\coolingconst]}$}{$\optimalbatcht{\iterstep+1} \gets  \optimalbatchtp{\iterstep+1}$} \label{inalg:worsestate2}
\Else{\label{inalg:staysame}$\optimalbatcht{\iterstep+1} \gets  \optimalbatcht{\iterstep}$}  \label{inalg:staysame2}}
} 
 \Return{$\optimalbatcht{\iterstep^*}$ \text{with the lowest } $\statesa{\iterstep^*}$.}
\end{algorithm}

We briefly describe the steps in Algorithm \ref{alg:simannealing}.
Line \ref{inalg:propose} finds a proposed \productwheel \ $\optimalbatchtp{\iterstep+1}$ given a current product wheel $\optimalbatcht{\iterstep}.$ For example, we can randomly perturb each $\bat{i}(k)$ to arrive at $\optimalbatchtp{\iterstep+1}$.  Lines \ref{inalg:check} - \ref{inalg:check2} then simulate $\numslot$ periods of production with demand schedule $\alldemand$ and the product wheel $\optimalbatchtp{\iterstep+1}$ under the skipping policy $\policy.$  Because $\pi$ is non-anticipatory, all skipping coefficients are determined by product wheel $\optimalbatchtp{\iterstep+1}$ during this simulated production run.  If, during this simulation, the \costconstraint \ or \inventoryconstraint \ are violated, a new product wheel is again proposed.

 In Lines \ref{inalg:betterstate} - \ref{inalg:betterstate2}, if the proposed product wheel has a shorter \avgwheeltime, we take the proposed product wheel as the next product wheel. However, as shown in lines \ref{inalg:worsestate} - \ref{inalg:worsestate2}, if the proposed product wheel has a longer \avgwheeltime, we accept the proposed product wheel with the probability of $\exp{[- \prtr{\statesap{\iterstep+1} -\statesa{\iterstep} }/\coolingconst]}$.  We otherwise keep the current product wheel, as shown in Lines \ref{inalg:staysame} - \ref{inalg:staysame2}.
 
 The entire process described above is repeated $\tmaxsa$ times, and all feasible product wheels found during these iterations are recorded.  The \productwheel \ that yields the smallest \avgwheeltime \ is returned.
 
    \section{Numerical Results}\label{sec:num}

In this section, we compare the performance of the optimization approaches of Section III using synthetic data.  In particular, we investigate how \costtolerance s, \changeoverc s, and \ \inventoryc s \  affect the \avgwheeltime s.  As shorthand, we refer to the integer linear programming strategy for optimization with no item skipping as the ILP approach, and the simulated annealing strategy for optimization with non-anticipatory skipping \eqref{eq:trigger} as the SA approach.

%%%%%%%%%%%%%%%%%%%%
\subsection{Synthetic data}

To generate the synthetic data, a demand schedule for three items over 24 time periods is sampled from a normal distribution with parameters given in Table \ref{table:syn data}.  To test robustness, five such demand schedules are generated.  Values for other parameters are also shown in Table \ref{table:syn data}.  For all simulations $\timeslot = 400$, and for the comparisons of \changeoverc \ and \inventoryc, we fix a cost tolerance of $\ct =40000$.  These synthetic values are chosen to be comparable to some particular manufacturing processes but have no special significance beyond illustration of trends and comparison of ILP and SA approaches.

%%%%%%%%%%%%%%%%%%%%%%%%
\subsection{Cost tolerance}
An illustration of the effect of cost tolerance on \avgwheeltime \ is shown in Figures \ref{fig:ilp-1} and \ref{fig:sa-1}.  Optimal \avgwheeltime s are found (when feasible) across 10 values of  $\ct$.  For low cost tolerance, the ILP approach finds no feasible solutions (represented by $\times$ in the figure and jittered for visibility).  The SA approach succeeds in finding optimal solutions with high \avgwheeltime s for low cost tolerance $\ct$; this is possible in part due to the lower \changeoverc \ of a large product wheel.  For both approaches, \avgwheeltime s\ clearly decrease as the cost tolerance goes up.  In particular, allowing for greater total \changeoverc s\ makes leaner production possible. Overall, the SA approach results in significantly lower \avgwheeltime \ than the ILP approach at comparable cost tolerances.

%%%%%%%%%%%%%%%%%%%
\subsection{Set-up cost}

An illustration of the effect of \changeoverc \ on \avgwheeltime \ is shown in Figures \ref{fig:ilp-2} and \ref{fig:sa-2}.  Optimal \avgwheeltime s are found (when feasible) across 10 multiples of \changeoverc \  $\ccost{i}$ in Table \ref{table:syn data}, all at a fixed cost tolerance of $\ct = 40000$.  As expected, \avgwheeltime \ increases as \changeoverc \ increases.   The SA approach finds feasible solutions for larger values of \changeoverc \ than the ILP approach, but with a sharp increase in \avgwheeltime \ in some cases.  Again, the SA approach results in significantly lower \avgwheeltime \ than the ILP approach at comparable \changeoverc.

%%%%%%%%%%%%%%%%%%%%%%
\subsection{Inventory cost}

An illustration of the effect of \inventoryc \  on \avgwheeltime \ is shown in Figures \ref{fig:ilp-3} and \ref{fig:sa-3}.  Similar to the previous examples, optimal \avgwheeltime s are found across 10 multiples of \inventoryc s\  $\icost{i}$ in Table \ref{table:syn data}, all at a fixed cost tolerance of $\ct = 40000$.  Increasing \inventoryc \ does not produce as strong an increase in \avgwheeltime \ as seen in the previous illustrations.  As expected, the SA approach again results in significantly lower \avgwheeltime \ than the ILP approach at comparable \inventoryc.

\begin{table*}%[]
\centering
\begin{tabular}{|c|c|c|c|c|c|c|c|c|}
\hline
Items       & Mean Monthly Demand & Std Dev Demand & $\ccost{i}$  & $\icost{i}$  & $\timebat{i}$ & $\unitbat{i}$ & $\invent{i}{0}$ & $\tp{i}$  \\ \hline
Item \#1 &	900 &	100 &	100 &	0.1	 & 1 &	10	& 0	& 100 \\ \hline
Item \#2 &	700	 & 140	& 100	& 0.15	& 1	& 15	& 0	&100 \\ \hline
Item \#3 &	1000	& 200	& 100	& 0.1	& 2	& 15	&0	&100 \\ \hline
\end{tabular}
\caption{Parameters for generating synthetic data.}
\label{table:syn data}
\end{table*}

\begin{figure*}
 \centering
     \begin{subfigure}[b]{0.42\linewidth}
 
         \includegraphics[width=\linewidth]{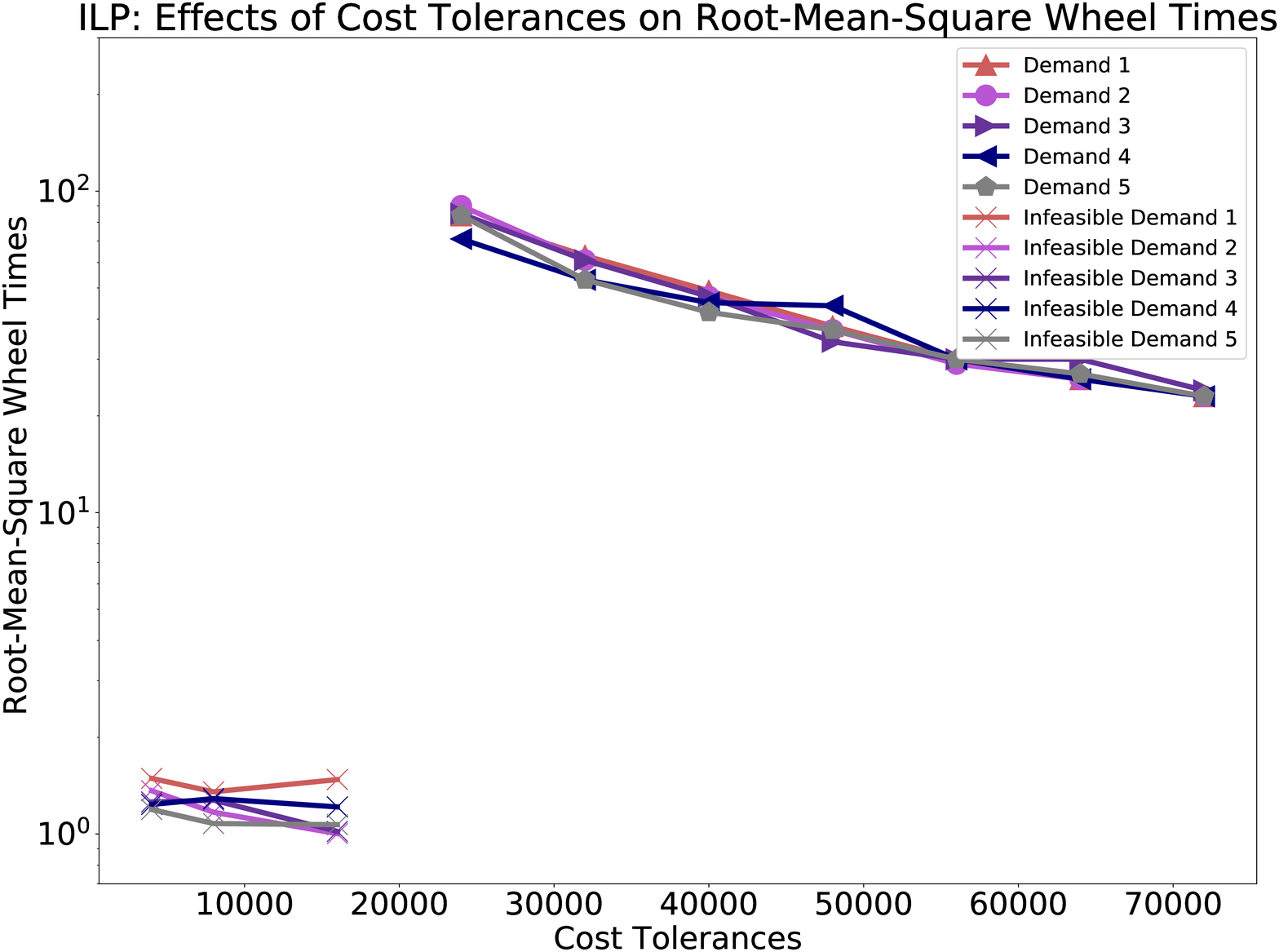}
         \caption{ILP: Effects of Cost Tolerances on Root-Mean-Square Wheel Times}
         \label{fig:ilp-1}
     \end{subfigure}
\hfill
       \begin{subfigure}[b]{0.42\linewidth}
     
         \includegraphics[width=\linewidth]{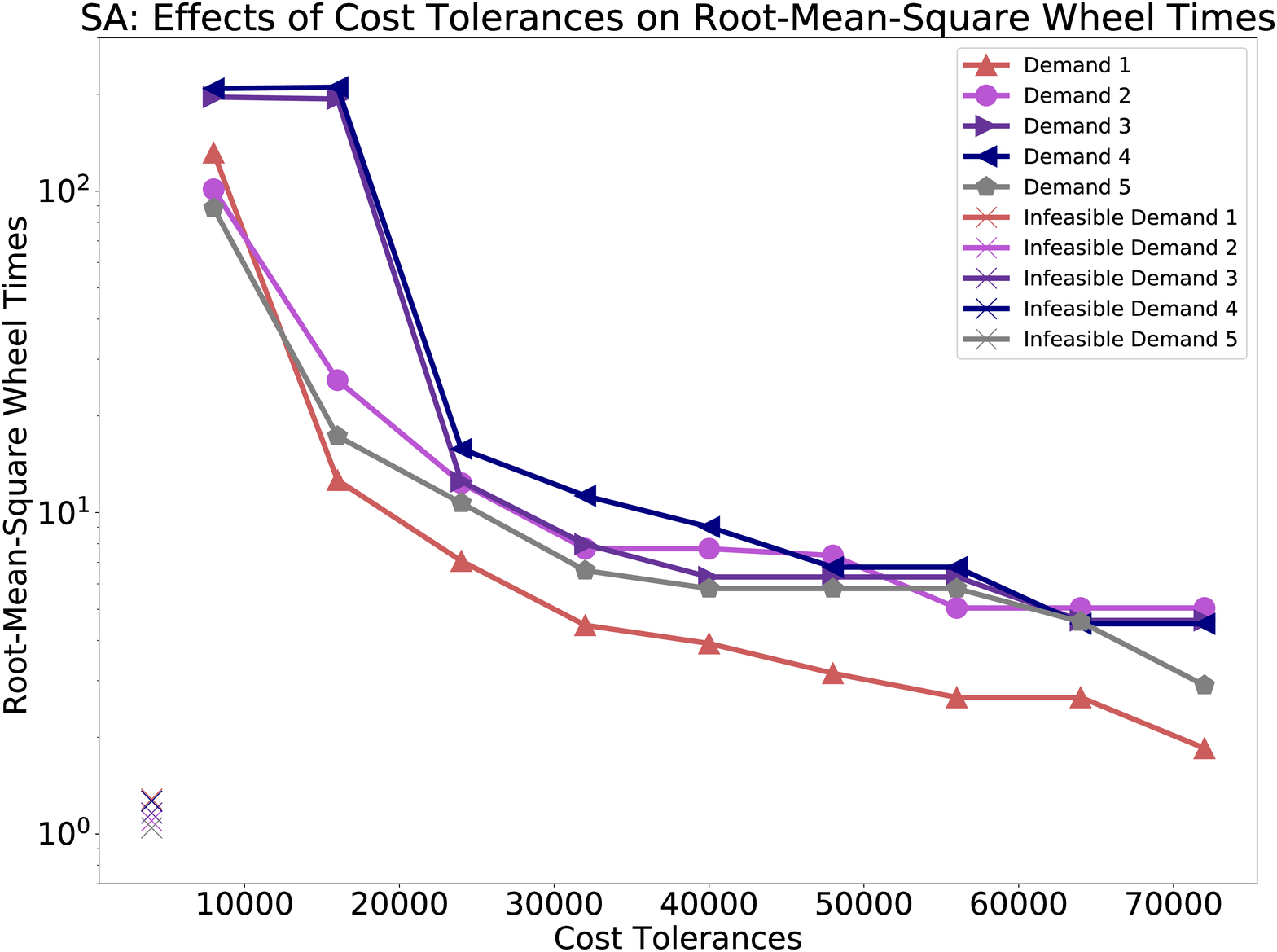}
         \caption{SA: Effects of Cost Tolerances on Root-Mean-Square Wheel Times}
         \label{fig:sa-1}
     \end{subfigure}
     
     \begin{subfigure}[b]{0.42\linewidth}
   
         \includegraphics[width=\linewidth]{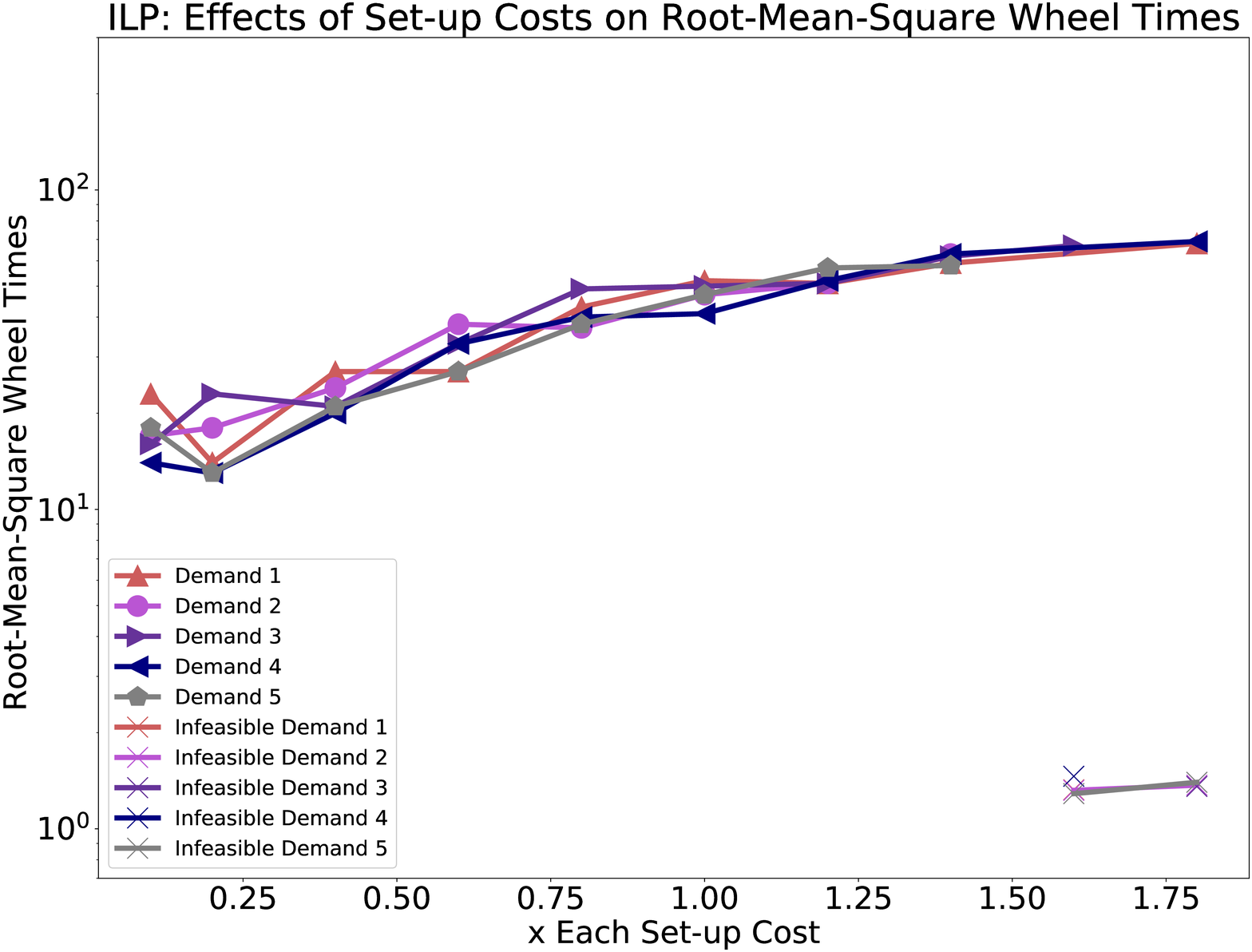}
         \caption{ILP: Effects of Set-up Costs on Root-Mean-Square Wheel Times}
         \label{fig:ilp-2}
     \end{subfigure}
 \hfill
  \begin{subfigure}[b]{0.42\linewidth}

         \includegraphics[width=\linewidth]{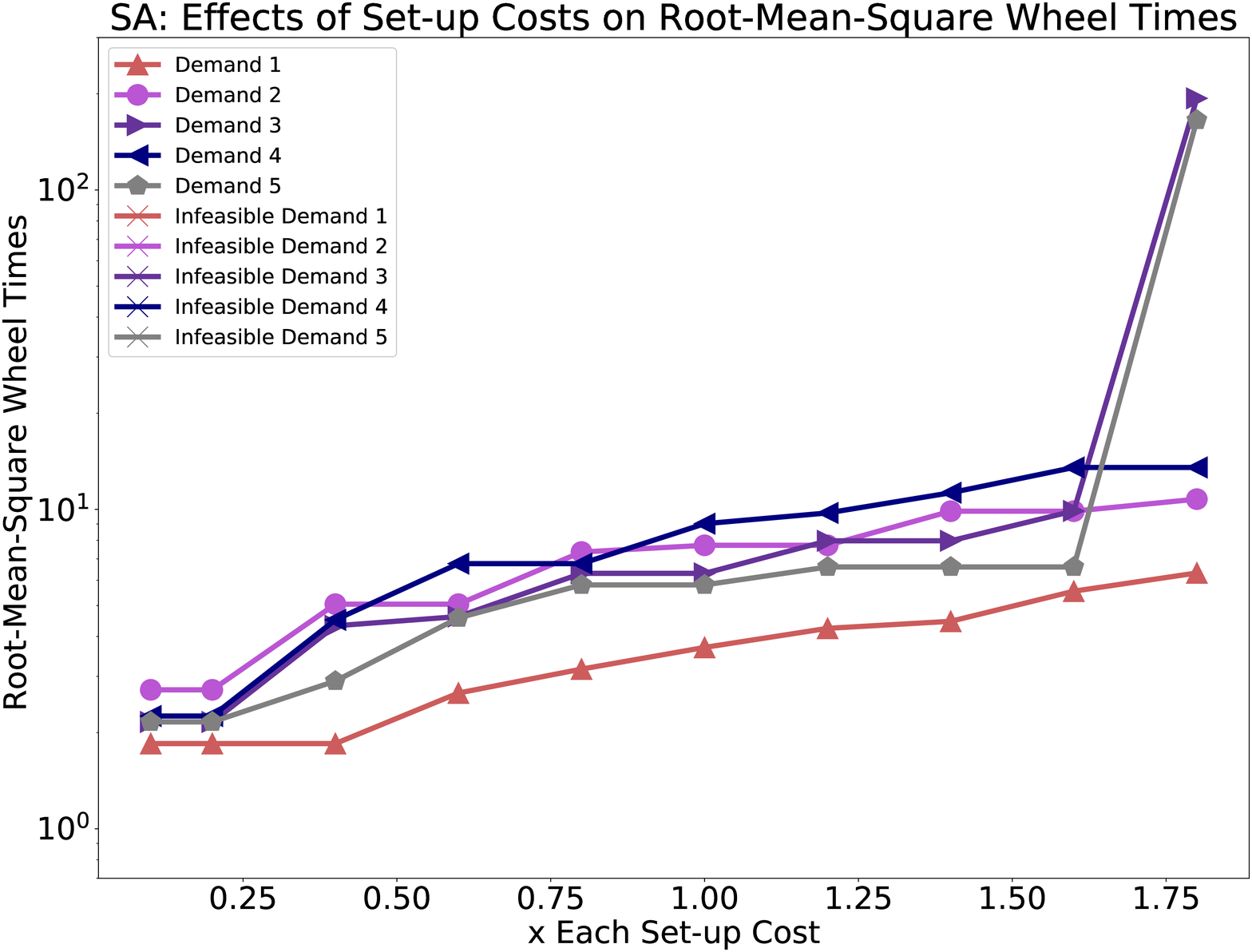}
         \caption{SA: Effects of Set-up Costs on Root-Mean-Square Wheel Times}
         \label{fig:sa-2}
     \end{subfigure}
     
\begin{subfigure}[b]{0.42\linewidth}
    
         \includegraphics[width=\linewidth]{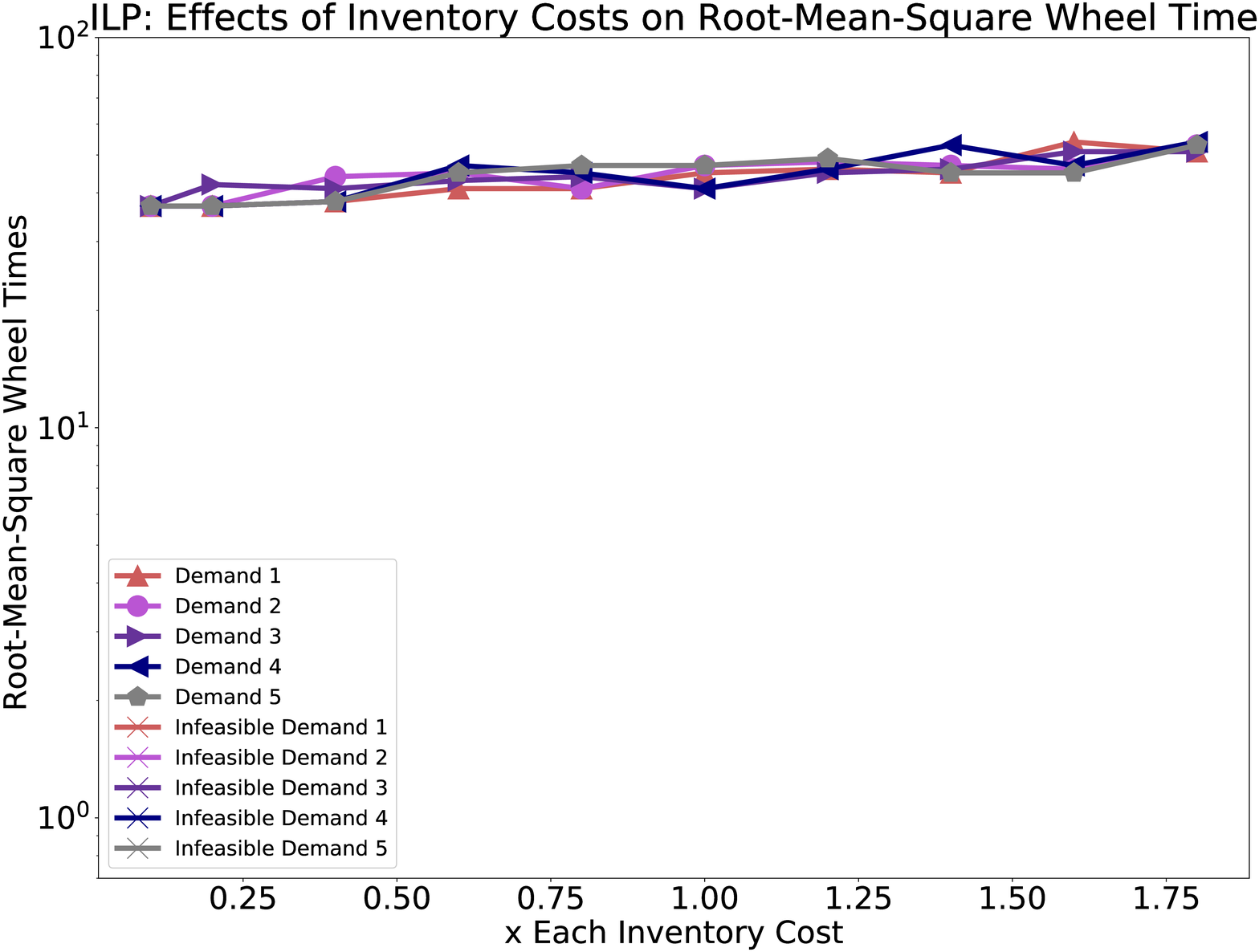}
         \caption{ILP: Effects of Inventory Costs on Root-Mean-Square Wheel Time}
         \label{fig:ilp-3}
     \end{subfigure}
  \hfill
     \begin{subfigure}[b]{0.42\linewidth}

         \includegraphics[width=\linewidth]{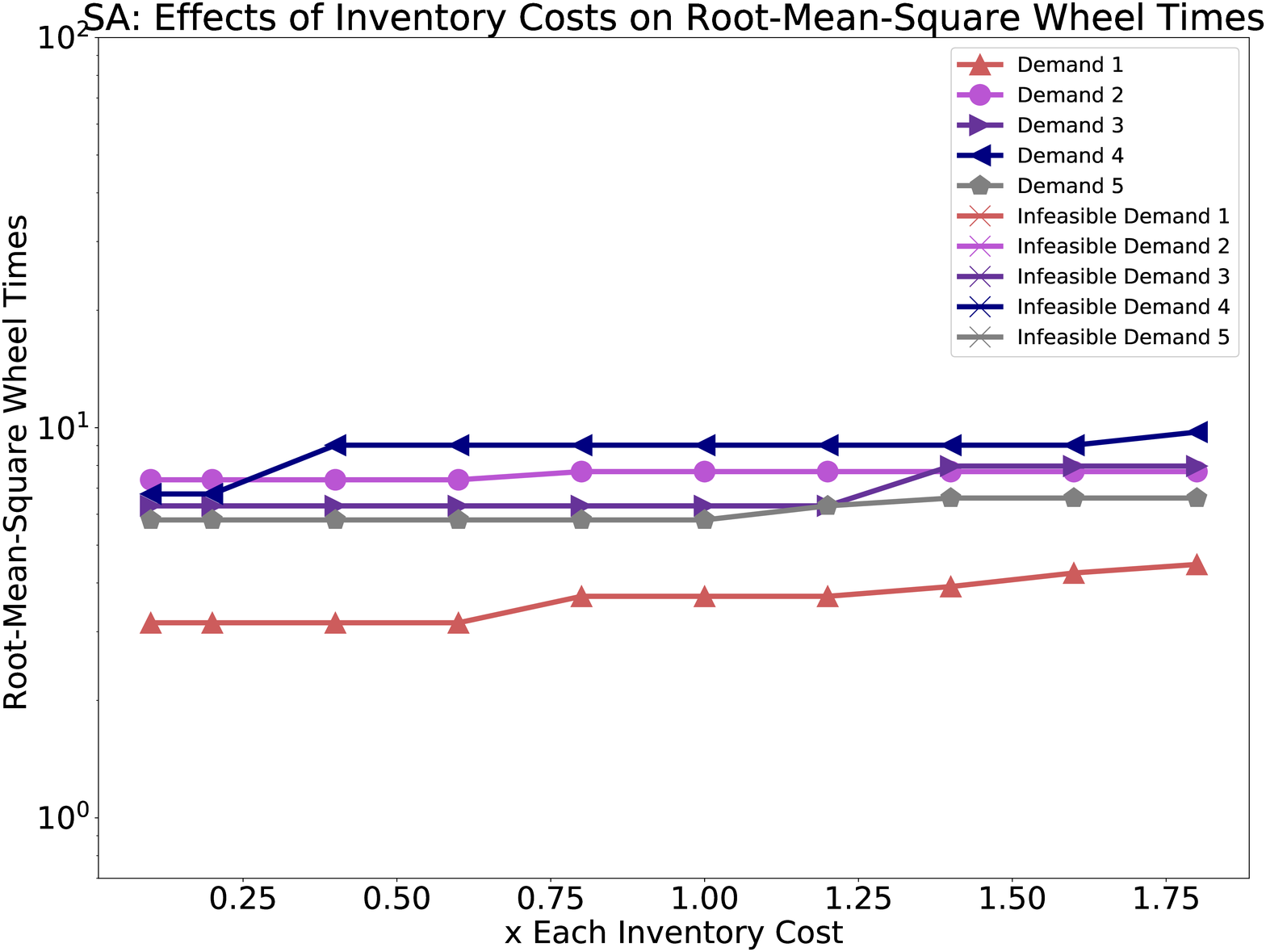}
         \caption{SA: Effects of Inventory Costs on Root-Mean-Square Wheel Times}
         \label{fig:sa-3}
     \end{subfigure}
        \caption{Optimal wheel time from SA and ILP approaches under varying \costtolerance, \changeoverc, and \inventoryc.}
        \label{fig:numerical result}
\end{figure*}
    \section{Conclusions And Future Research}
 \label{sec:conclusion}
 
In this work, we introduced a mathematical framework for the product wheel methodology, and cast the problem of creating a lean manufacturing schedule as a nonconvex optimization problem.  We applied simplifying assumptions to the skipping policy for the product wheel and created reformulations of this optimization problem solvable via integer linear programming and simulated annealing.  Finally, we evaluated the performance of these optimization strategies on a synthetic data set, investigating how different costs affect leanness.

Several future research directions are apparent.  (i) The simplifying assumptions of no skipping and non-anticipatory skipping were essential for the optimization strategies presented.  It would be valuable to explore strategies for optimization outside those assumptions.  In particular, a global optimization over all possible skipping choices could achieve still shorter \wheeltime s.  (ii) Throughout this work, we have assumed full knowledge of demand schedules, which is unrealistic in most manufacturing scenarios.  Addressing the propagation of uncertainty through the optimization process, e.g., via Monte Carlo simulation or robust optimization, \cite{powell2019unified}, \cite{ben2009robust}, would have significant practical value.  (iii) Finally, our simulations make the overall advantage of the SA approach over the ILP approach apparent but do not explore more extreme demand schedules or costs (e.g., demand with large fluctuations or large disparities in costs between different items).  It would be worthwhile to do a more thorough numerical exploration of  these extremes to understand ideal cases to employ non-anticipative skipping.

     \section{Acknowledgements}
We are grateful to Joseph G. Caradimitropoulo, Luanqi Chen,
William Thompson,  Yangqiaoyu Zhou, Marcella Manivel and Libby Nachreiner for essential contributions to this project, including work on the problem formulation, linear programming based solutions, and skipping coefficient simulations.  Our correspondents at 3M, Mike Gerlach, Digital Manufacturing Technologies Manager; Jordan Reed, Master Production Planner; and Michael Muilenburg, Operational Technologies Manager, introduced us to the foundational ideas for this project and offered valuable feedback on early work.  Partial funding for the first author was provided by the Carleton College Towsley Endowment.

    % \bibliography{reference}
    % \bibliographystyle{IEEEtran}
    \bibliography{reference}
    \bibliographystyle{plain}

\end{document}